\documentclass[10pt]{article}
\usepackage{amsmath}
\usepackage{amssymb}
\usepackage{mathrsfs}
\usepackage{amsfonts}
\usepackage{mathrsfs,amscd,amssymb,amsthm,amsmath,bm,graphicx}
\usepackage{graphicx}
\usepackage{cite}
\usepackage{color}
\usepackage{enumitem}
\usepackage{setspace}
\makeatletter

\renewcommand{\@seccntformat}[1]{{\csname the#1\endcsname}{\normalsize.}\hspace{.5em}}
\makeatother

\def \[{\begin{equation}}
\def \]{\end{equation}}

\newtheorem{thm}{Theorem}[section]

\newtheorem{lem}[thm]{Lemma}

\begin{document}
\setlength{\baselineskip}{13pt}
\begin{center}{\Large \bf  Resistance distance-based graph invariants and the number of spanning trees
of linear crossed octagonal graphs
}

\vspace{4mm}

{\large Jing Zhao$^1$,~Jia-Bao Liu$^{1,2,*}$,~Sakander Hayat$^3$}\vspace{2mm}

{\small $^1$ School of Mathematics and Physics, Anhui Jianzhu
University, Hefei
230601, P.R. China\\
$^2$ School of Mathematics, Southeast University,
Nanjing 210096, P.R. China\\
$^3$ Faculty of Engineering Sciences, GIK Institute of Engineering Sciences and Technology, Topi 23460, Pakistan}
\vspace{2mm}
\end{center}

\footnotetext{E-mail address: zhaojing94823@163.com, liujiabaoad@163.com, sakander1566@gmail.com.}

\footnotetext{* Corresponding author.}

 {\noindent{\bf Abstract.}\ \  Resistance distance is a novel distance function,
 also a new intrinsic graph metric, which makes some extensions of ordinary distance.
 Let $O_n$ be a linear crossed octagonal graph. Recently, Pan and Li (2018) derived the closed formulas for the Kirchhoff index, multiplicative degree-Kirchhoff index and the number of spanning trees of $H_n$. They pointed that it is interesting to give the explicit formulas for the Kirchhoff and  multiplicative degree-Kirchhoff indices of $O_n$. Inspired by these, in this paper, two resistance
 distance-based graph invariants, namely, Kirchhoff and multiplicative
 degree-Kirchhoff indices are studied. We firstly determine formulas for
 the Laplacian (normalized Laplacian, resp.) spectrum of $O_n$.
Further, the formulas for those two resistance distance-based graph invariants and
spanning trees are given. More surprising, we find that the Kirchhoff (multiplicative degree-Kirchhoff, resp.) index
is almost one quarter to Wiener (Gutman, resp.) index of a linear crossed octagonal graph.

\noindent{\bf Keywords}:  Laplacian, Normalized Laplacian, Kirchhoff index, Multiplicative degree-Kirchhoff index, Spanning tree\vspace{2mm}

\noindent{\bf AMS subject classification:} 05C50,\ 05C90}

\section{Introduction} \ \ \ \ \ We will recall some background and notations that lead to the main
results in subsection 1.1. Besides, we present the main results of linear crossed
octagonal graph in the subsections 1.2 and 1.3.
\subsection{ Background}
\ \ \ \ \ Many practical problems in real world can be described
or characterized by graphs. In this paper, simple and undirected graphs are considered.
Assume that $G=(V_{G},E_{G})$ is a graph. $V_{G}=\{v_1, v_2,\ldots,v_n\}$ is the vertex-set
and $E_{G}=\{e_1,e_2,\ldots,e_m\}$ the edge-set of $G$. Further, $|V_G|$ and $|E_G|$
 are the order
and size of $G$, respectively. One follows \cite{F.R} for
terminology and notation that we omit here.

In graph theory, some parameters were used to
characterize the structure properties of graphs.
While many of them have practical applications in the fields
of physico-chemical, biological and etc., some only have pure
mathematics properties. Distance is one of the most studied
quantity in graph, such as the sum of the distance,
denoted by $W(G)=\sum_{i<j}d_{ij}$, was named Wiener index\cite{Wiener},
where $d_{ij}$ is the length of a shortest path between vertices $i$ and $j.$
In 1994, Gutman introduced the Gutman index\cite{I.G} as $Gut(G)=\sum_{i<j}d_id_jd_{ij}.$

The resistance distance\cite{KR}, a novel distance function, also a new
intrinsic graph metric, was proposed by Klein and Randi\'c in 1993 on the basis of the theory of electrical networks.
Denote $r_{ij}$ the resistance distance between any two vertices $i$ and $j$ in graphs by a unit resistor instead of each edge.
Klein and Randi\'{c} \cite{D.J.,D.} defined the Kirchhoff index on the basis of the resistance distance and given it by
$Kf(G)=\sum_{i<j}r_{ij}.$ The homologous consideration with Gutman index, Chen and Zhang~\cite{H.Y.} presented the multiplicative
degree-Kirchhoff index as $Kf^*(G)=\sum_{i<j}d_id_jr_{ij}.$ Indeed, the Kirchhoff and Wiener index coincide (multiplicative
degree-Kirchhoff and Gutman index, resp.), if $G$ is a tree.

The number of spanning trees\cite{F.R} of a given graph $G$ is the
number of subgraphs which contains each vertices of $G$ and all those subgraphs shall be trees.

In spectral graph theory, the interplays between the
structure properties and eigenvalues of $G$ were highly concerned. The adjacency matrix $A(G)=(a_{ij})_{n\times n}$ of $G$ is an $(0,1)$-matrix, $a_{ij}=1$ or $0$, if there exits an edge between vertices $i$ and $j$ or not. Define $d_G(v_i)$ ($d_i$ for simply) the degree of vertex $v_i$ and $D(G)=diag(d_1,d_2,\ldots,d_n)$ the degree matrix of
$G$. The
Laplacian matrix was defined as
\begin{eqnarray*}
\big(L(G)\big)_{ij}=
\begin{cases}
d_i,       & i=j; \\
-1,   & i\neq j,\ \ v_i\ \ and\ \ v_j\ \ are\ \ adjacency; \\
0,           &otherwise.
\end{cases}
\end{eqnarray*}

The normalized Laplacian
was given by
\begin{eqnarray*}
\big(\mathcal{L}(G)\big)_{ij}=
\begin{cases}
1,       & i=j; \\
-\frac{1}{\sqrt{d_id_j}},   & i\neq j,\ \ v_i\ \ and\ \ v_j\ \ are\ \ adjacency; \\
0,           &otherwise.
\end{cases}
\end{eqnarray*}

According to the definitions of Lpalcian and normalized Laplacian, one immediately gets
$L(G)=D(G)-A(G)$ and $\mathcal{L}(G)=D(G)^{-\frac{1}{2}}L(G)D(G)^{-\frac{1}{2}},$ see\cite{H.H,J.L,J,H.,Mohar} for details.

Gutman et al.\cite{Gut} explored that the Kirchhoff and quasi-Wiener index coincide, that is
\begin{eqnarray}
Kf(G)&=&n\sum_{k=2}^n\frac{1}{\mu_k},
\end{eqnarray}
where $0=\mu_1<\mu_2\leq\cdots \leq \mu_n$ are the Laplacian eigenvalues of
$L(G)$.

Based on the normalized Laplacian, Chen and Zhang\cite{H.Y.} in 2007 proved that
\begin{eqnarray}
Kf^*(G)&=&2m\sum_{k=1}^n\frac{1}{\lambda_k},
\end{eqnarray}
where $\lambda_1\leq\lambda_2\leq\cdots\leq\lambda_n$ are the normalized Laplacian
eigenvalues of $\mathcal{L}(G)$.

Since the extensively applications of the Kirchhoff index in the fields of physics, chemistry and networks science, such as,
the smaller the Kirchhoff index, the better the connectivity of a given graph. It is more meaningful to determine the Kirchhoff or multiplicative degree-Kirchhoff index. Up to date, the Kirchhoff or multiplicative degree-Kirchhoff index have already determined of some
classes of graphs, see\cite{J2016,Y2017,chun2018,Zhu2018,S2018,X2019,Zhao}. In this paper, we only discuss a class of linear graphs. In 2007, Yang et al.\cite{Y2007} obtained the Laplacian spectrum and Kirchhoff index of linear hexagonal graph. The normalized Laplacian and multiplicative degree-Kirchhoff index of linear hexagonal graph were determined by Huang et al.\cite{HL}. Pan et al.\cite{P2018} constructed a linear crossed hexagonal graph, the Kirchhoff and  multiplicative degree-Kirchhoff indices were explored. Meanwhile, they pointed out that it is interesting to give the explicit formulas for the Kirchhoff and  multiplicative degree-Kirchhoff indices of linear crossed octagonal graph, see Figure 1.

Inspired by that interesting problem, we will explore the explicit formulas for the Kirchhoff and  multiplicative degree-Kirchhoff indices of linear crossed octagonal graph.

\begin{figure}[htbp]
\centering\includegraphics[width=12cm,height=3cm]{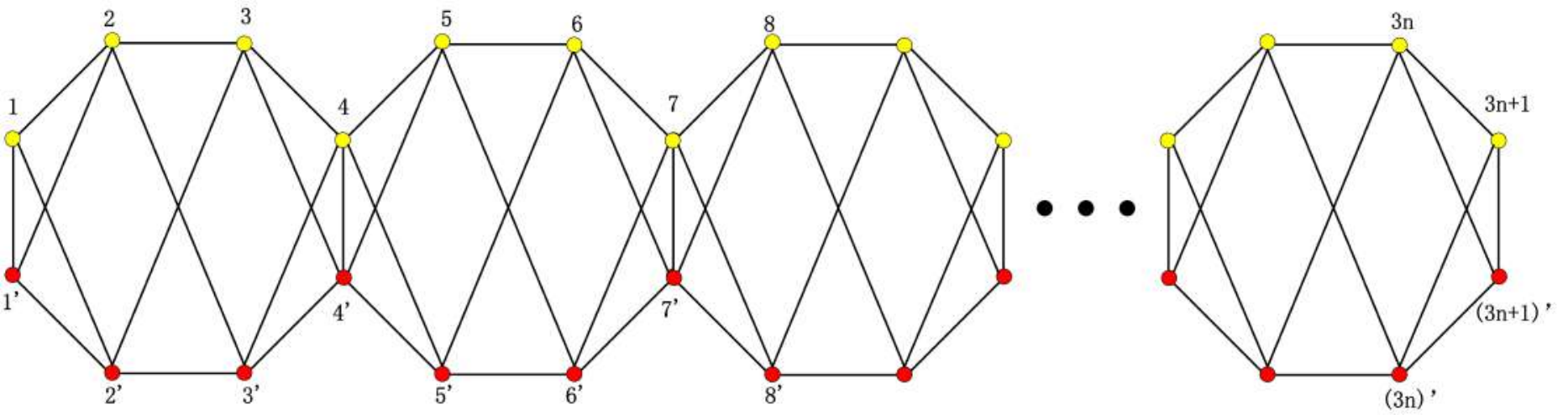}
\caption{Linear crossed octagonal graph.}\label{1}
\end{figure}

The rest of the paper are concluded as follows. In this context, a method  that has a big effect on the proofs of main results is presented in section 2. The proof of the Kirchhoff index of linear crossed octagonal graph is determined on the basis of the Laplacian spectrum. Also, one proves that the ratio between the Kirchhoff and Wiener index is almost one quarter of linear crossed octagonal graph in section 3. One determines the multiplicative degree-Kirchhoff index, the number of spanning trees and the ratio between the multiplicative degree-Kirchhoff and Gutman index is almost one quarter of linear crossed octagonal graph in section 4. The conclusion of whole paper is concluded in section 5.

\subsection{ Kirchhoff index of linear crossed octagonal graph}
\ \ \ \ \ In this subsection, one presents the Kirchhoff index and the ratio  between the Kirchhoff and Wiener index of linear crossed octagonal graph as follows.
\begin{thm} Suppose that $O_n$ are the linear crossed octagonal graphs. Then
\begin{eqnarray*}
Kf(O_n)=\frac{27n^3+51n^2+26n+4}{6}.
\end{eqnarray*}
\end{thm}

To our surprise, one finds the Kirchhoff index is almost one quarter to its Wiener index of linear crossed octagonal graphs.

\begin{thm} Suppose that $O_n$ are linear crossed octagonal graphs. Then
\begin{eqnarray*}
\lim_{n\rightarrow \infty}\frac{Kf(O_n)}{W(O_n)}=\frac{1}{4}.
\end{eqnarray*}
\end{thm}

\subsection{ Multiplicative degree-Kirchhoff index of linear crossed octagonal graph}
\ \ \ \ \ In this subsection, one proposes the multiplicative degree-Kirchhoff index, the number of spanning trees and the ratio between the multiplicative degree-Kirchhoff and Gutman index of linear crossed octagonal graph as below.

\begin{thm} For linear crossed
octagonal graphs $O_n$,
\begin{eqnarray*}
Kf^*(O_n)=\frac{507n^3+559n^2+204n+8}{6}.
\end{eqnarray*}
\end{thm}

At this place, the number of spanning trees of linear crossed octagonal graphs is given as below.

\begin{thm} For linear crossed
octagonal graphs $O_n,$
\begin{eqnarray*}
\tau(O_n)=2^{8n+2}\cdot3^{n-1}.
\end{eqnarray*}
\end{thm}

The same to the relationship between the Kirchhoff and Wiener index, the multiplicative degree-Kirchhoff is almost one quarter to its Gutman index of linear crossed octagonal graphs.

\begin{thm} Suppose that $O_n$ are linear crossed octagonal graphs. Then
\begin{eqnarray*}
\lim_{n\rightarrow \infty}\frac{Kf^*(O_n)}{Gut(O_n)}=\frac{1}{4}.
\end{eqnarray*}
\end{thm}

\section{Preliminary}
\ \ \ \ In this section, we will recall an important method which was presented in \cite{Yang}.
Based on this method, some special matrices of linear crossed octagonal graphs will be given.

Let $P_{R}(x)=\det (xI-R)$ be the characteristic
polynomial of the $n\times n$ matrix $R$. Assume that $O_n$ is the linear crossed octagonal graph, of which is generated by adding three pairs of crossed edges of each octagonal, see Figure 1. One can straightforward to check that $|V(G)|=6n+2$ and $|E(G)|=13n+1.$ Assume that
$V_1=\{1,2,\ldots,3n+1\}, V_2=\{1',2',\ldots,(3n+1)'\}$. We immediately obtain that $\pi=(1,1')(2,2')\cdots(3n+1,(3n+1)')$ is an automorphism.  Hence,
the Laplacian and normalized Laplacian matrices of $O_n$ can be expressed by
 block matrices as below.
\begin{equation*}
L(O_n)=\left(
  \begin{array}{cc}
  L_{V_1V_1}& L_{V_1V_2}\\
 L_{V_2V_1}& L_{V_2V_2}\\
  \end{array}
\right),~
\mathcal{L}(O_n)=\left(
  \begin{array}{cc}
  \mathcal{L}_{V_1V_1}& \mathcal{L}_{V_1V_2}\\
 \mathcal{L}_{V_2V_1}& \mathcal{L}_{V_2V_2}\\
  \end{array}
\right).
\end{equation*}
One easily obtains that $L_{V_1V_1}=L_{V_2V_2}$,
$L_{V_1V_2}=L_{V_2V_1}$,
$\mathcal{L}_{V_1V_1}=\mathcal{L}_{V_2V_2}$ and
$\mathcal{L}_{V_1V_2}=\mathcal{L}_{V_2V_1}$.

Suppose that
\begin{equation*}
P=\left(
  \begin{array}{cc}
  \frac{1}{\sqrt{2}}I_{n}& \frac{1}{\sqrt{2}}I_{n}\\
  \frac{1}{\sqrt{2}}I_{n}& -\frac{1}{\sqrt{2}}I_{n}
  \end{array}
\right).
\end{equation*}
Consequently
\begin{equation*}
PL(O_n)P'=\left(
  \begin{array}{cc}
  L_A& 0\\
    0& L_{s}
  \end{array}
\right),~
P\mathcal{L}(O_n)P'=\left(
  \begin{array}{cc}
  \mathcal{L}_A& 0\\
    0& \mathcal{L}_{s}
  \end{array}
\right),
\end{equation*}
where $P'$ is the transposition of $P$.

Further,
\begin{eqnarray}
L_A=L_{V_1V_1}+L_{V_1V_2},
~L_s=L_{V_1V_1}-L_{V_1V_2}.\\
\mathcal{L}_A=\mathcal{L}_{V_1V_1}+\mathcal{L}_{V_1V_2},
~\mathcal{L}_s=\mathcal{L}_{V_1V_1}-\mathcal{L}_{V_1V_2}.
 \end{eqnarray}

According to (2.3), one gets block matrices of $L_{V_1V_1}$ and $L_{V_1V_2}$ in term of the definition of Laplacian as follows.
\begin{eqnarray*}
 L_{V_1V_1}&=&(l_{ij})_{(3n+1)\times(3n+1)}\\&=&
\left(
  \begin{array}{cccccccccc}
3 & -1 & 0& 0& 0& \cdots& 0& 0& 0& 0\\
-1 & 4 & -1& 0& 0& \cdots& 0& 0& 0& 0\\
0& -1 & 4 & -1 & 0& \cdots& 0& 0& 0& 0\\
0& 0 & -1 & 5 &-1 &  \cdots& 0& 0& 0& 0\\
0& 0 & 0 &-1 & 4 &\cdots  & 0& 0& 0& 0\\
\vdots& \vdots& \vdots& \vdots& \vdots&\ddots  & \vdots& \vdots& \vdots& \vdots\\
0& 0& 0& 0& 0& \cdots & 5& -1 & 0& 0\\
0& 0& 0& 0& 0& \cdots& -1 &4 &-1 & 0\\
0& 0& 0& 0& 0& \cdots& 0&-1 &4 &-1  \\
0& 0& 0& 0& 0& \cdots& 0& 0&-1 &3  \\
  \end{array}
\right)_{(3n+1)\times (3n+1)},
\end{eqnarray*}
where $l_{11}=l_{3n+1,3n+1}=3$,
$l_{3i-1,3i-1}=l_{3i,3i}=4$ for $i\in\{1,2,\ldots,n\}$ and
$l_{3i+1,3i+1}=5$ for $i\in\{1,2,\ldots,n-1\}$. Also,
$l_{i,i+1}=l_{i+1,i}=-1$ for $i\in\{1,2,\ldots,3n\}$.

\begin{eqnarray*}
 L_{V_1V_2}&=&(l_{ij})_{(3n+1)\times(3n+1)}\\&=&
\left(
  \begin{array}{cccccccccc}
-1 & -1 & 0& 0& 0& \cdots& 0& 0& 0& 0\\
-1 & 0 & -1& 0& 0& \cdots& 0& 0& 0& 0\\
0& -1 & 0 & -1 & 0& \cdots& 0& 0& 0& 0\\
0& 0 & -1 & -1 &-1 &  \cdots& 0& 0& 0& 0\\
0& 0 & 0 &-1 & 0 &\cdots  & 0& 0& 0& 0\\
\vdots& \vdots& \vdots& \vdots& \vdots&\ddots  & \vdots& \vdots& \vdots& \vdots\\
0& 0& 0& 0& 0& \cdots & -1& -1 & 0& 0\\
0& 0& 0& 0& 0& \cdots& -1 &0 &-1 & 0\\
0& 0& 0& 0& 0& \cdots& 0&-1 &0 &-1  \\
0& 0& 0& 0& 0& \cdots& 0& 0&-1 &-1  \\
  \end{array}
\right)_{(3n+1)\times (3n+1)},
\end{eqnarray*}
where $l_{11}=l_{3n+1,3n+1}=-1$,
$l_{3i-1,3i-1}=l_{3i,3i}=0$ for $i\in\{1,2,\ldots,n\}$ and
$l_{3i+1,3i+1}=-1$ for $i\in\{1,2,\ldots,n-1\}$. Also,
$l_{i,i+1}=l_{i+1,i}=-1$ for $i\in\{1,2,\ldots,3n\}$.

For the normalized Laplacian of linear crossed octagonal graphs, the blocks matrices
of $\mathcal {L}_{V_1V_1}$ and $\mathcal {L}_{V_1V_2}$ are as below.
\begin{eqnarray*}
 \mathcal {L}_{V_1V_1}&=&(l_{ij})_{(3n+1)\times(3n+1)}\\&=&
\left(
  \begin{array}{cccccccccc}
1 & -\frac{1}{\sqrt{12}} & 0& 0& 0& \cdots& 0& 0& 0& 0\\
-\frac{1}{\sqrt{12}} & 1 & -\frac{1}{4}& 0& 0& \cdots& 0& 0& 0& 0\\
0& -\frac{1}{4} & 1 & -\frac{1}{\sqrt{20}} & 0& \cdots& 0& 0& 0& 0\\
0& 0& -\frac{1}{\sqrt{20}}&1 &-\frac{1}{\sqrt{20}} &  \cdots& 0& 0& 0& 0\\
0& 0& 0&-\frac{1}{\sqrt{20}} & 1 &\cdots  & 0& 0& 0& 0\\
\vdots& \vdots& \vdots& \vdots& \vdots&\ddots  & \vdots& \vdots& \vdots& \vdots\\
0& 0& 0& 0& 0& \cdots & 1&-\frac{1}{\sqrt{20}} & 0& 0\\
0& 0& 0& 0& 0& \cdots& -\frac{1}{\sqrt{20}}&1 &-\frac{1}{4} & 0\\
0& 0& 0& 0& 0& \cdots& 0&-\frac{1}{4} &1 &-\frac{1}{\sqrt{12}}  \\
0& 0& 0& 0& 0& \cdots& 0& 0&-\frac{1}{\sqrt{12}} &1  \\
  \end{array}
\right)_{(3n+1)\times (3n+1)},
\end{eqnarray*}
where $l_{i,i}=1$, for $i\in\{1,2,\ldots,3n+1\}$. Also,
$l_{12}=l_{21}=l_{3n,3n+1}=l_{3n+1,3n}=-\frac{1}{\sqrt{12}}$,
$l_{3i,3i-1}=l_{3i-1,3i}=-\frac{1}{4}$, for $i\in\{1,2,\ldots,n\}$, $l_{3i+1,3i}=l_{3i,3i+1}=l_{3i+2,3i+1}=l_{3i+1,3i+2}=-\frac{1}{\sqrt{20}}$,
for $i\in\{1,2,\ldots,n-1\}$.

\begin{eqnarray*}
 \mathcal{L}_{V_1V_2}&=&(l_{ij})_{(3n+1)\times(3n+1)}\\&=&
\left(
  \begin{array}{cccccccccc}
-\frac{1}{3} & -\frac{1}{\sqrt{12}} & 0& 0& 0& \cdots& 0& 0& 0& 0\\
-\frac{1}{\sqrt{12}} & 0 & -\frac{1}{4}& 0& 0& \cdots& 0& 0& 0& 0\\
0& -\frac{1}{4} & 0 & -\frac{1}{\sqrt{20}} & 0& \cdots& 0& 0& 0& 0\\
0& 0& -\frac{1}{\sqrt{20}}&-\frac{1}{5} &-\frac{1}{\sqrt{20}} &  \cdots& 0& 0& 0& 0\\
0& 0& 0&-\frac{1}{\sqrt{20}} & 0 &\cdots  & 0& 0& 0& 0\\
\vdots& \vdots& \vdots& \vdots& \vdots&\ddots  & \vdots& \vdots& \vdots& \vdots\\
0& 0& 0& 0& 0& \cdots & -\frac{1}{5} &-\frac{1}{\sqrt{20}} & 0& 0\\
0& 0& 0& 0& 0& \cdots& -\frac{1}{\sqrt{20}}&0 &-\frac{1}{4} & 0\\
0& 0& 0& 0& 0& \cdots& 0&-\frac{1}{4} &0 &-\frac{1}{\sqrt{12}}  \\
0& 0& 0& 0& 0& \cdots& 0& 0&-\frac{1}{\sqrt{12}} &-\frac{1}{3}  \\
  \end{array}
\right)_{(3n+1)\times (3n+1)},
\end{eqnarray*}
where $l_{1,1}=l_{3n+1,3n+1}=-\frac{1}{3}$, $l_{3i-1,3i-1}=l_{3i,3i}=0$ for $i\in\{1,2,\ldots,n\}$ and
$l_{3i+1,3i+1}=-\frac{1}{5}$ for $i\in\{1,2,\ldots,n-1\}$. Also,
$l_{12}=l_{21}=l_{3n,3n+1}=l_{3n+1,3n}=-\frac{1}{\sqrt{12}}$,
$l_{3i,3i-1}=l_{3i-1,3i}=-\frac{1}{4}$, for $i\in\{1,2,\ldots,n\}$, $l_{3i+1,3i}=l_{3i,3i+1}=l_{3i+2,3i+1}=l_{3i+1,3i+2}=-\frac{1}{\sqrt{20}}$,
for $i\in\{1,2,\ldots,n-1\}$.

In what follows, the lemmas that we present are necessarily part
in the proofs of the main results.

\begin{lem}\cite{Y2007,HL} Assume that $L_A, L_S, \mathcal{L}_A$, and $\mathcal{L}_{S}$ are obtained as mentioned above. Then
$$P_{L(G)}(x)=P_{L_A}(x)P_{L_{s}}(x),~P_{\mathcal{L}(G)}(x)=P_{\mathcal{L}_A}(x)P_{\mathcal{L}_{s}}(x).$$

\end{lem}
\begin{lem}\cite{F.R}
Assume that $G$ is a graph with $|V_G|=n$ and $|E_G|=m$. Then
$$\prod_{i=1}^nd_G(v_i)\prod_{i=2}^n\lambda_i=2m\tau(G),$$
where
$\tau(G)$ is the number of spanning trees of $G$.
\end{lem}

\section{Proofs of Theorems 1.1 and 1.2}
\ \ \ \ In this section, we will first give the proof of Theorem 1.1. According to (1.1), one will devote to figure out the Laplacian spectrum of $O_n$ due to its crucial role in the proof of Theorem 1.1.
\vskip 0.3 cm

\noindent{\bf Proof of Theorem 1.1.}
In view of (2.3), one gets
\begin{eqnarray*}
L_{A}&=& \left(
  \begin{array}{cccccccccc}
2 & -2 & 0& 0& 0& \cdots& 0& 0& 0& 0\\
-2 & 4 & -2& 0& 0& \cdots& 0& 0& 0& 0\\
0& -2 & 4 & -2 & 0& \cdots& 0& 0& 0& 0\\
0& 0 & -2 & 4 &-2 &  \cdots& 0& 0& 0& 0\\
0& 0 & 0 &-2 & 4 &\cdots  & 0& 0& 0& 0\\
\vdots& \vdots& \vdots& \vdots& \vdots&\ddots  & \vdots& \vdots& \vdots& \vdots\\
0& 0& 0& 0& 0& \cdots & 4& -2 & 0& 0\\
0& 0& 0& 0& 0& \cdots& -2 &4 &-2 & 0\\
0& 0& 0& 0& 0& \cdots& 0&-2 &4 &-2  \\
0& 0& 0& 0& 0& \cdots& 0& 0&-2 &2  \\
  \end{array}
\right)_{(3n+1)\times (3n+1)},
\end{eqnarray*}
where $l_{11}=l_{3n+1,3n+1}=2$,
$l_{i,i}==4$ for $i\in\{2,3,\ldots,3n\}$. Also,
$l_{i,i+1}=l_{i+1,i}=-2$ for $i\in\{1,2,\ldots,3n\}$.

\begin{eqnarray*}
L_S&=& \left(
  \begin{array}{cccccccccc}
4 & 0 & 0& 0& 0& \cdots& 0& 0& 0& 0\\
0 & 4 & 0& 0& 0& \cdots& 0& 0& 0& 0\\
0& 0 & 4 & 0 & 0& \cdots& 0& 0& 0& 0\\
0& 0 & 0 & 6 &0 &  \cdots& 0& 0& 0& 0\\
0& 0 & 0 &0 & 4 &\cdots  & 0& 0& 0& 0\\
\vdots& \vdots& \vdots& \vdots& \vdots&\ddots  & \vdots& \vdots& \vdots& \vdots\\
0& 0& 0& 0& 0& \cdots & 6& 0 & 0& 0\\
0& 0& 0& 0& 0& \cdots& 0 &4 &0 & 0\\
0& 0& 0& 0& 0& \cdots& 0&0 &4 &0  \\
0& 0& 0& 0& 0& \cdots& 0& 0&0 &4  \\
  \end{array}
\right)_{(3n+1)\times (3n+1)},
\end{eqnarray*}
where $l_{11}=l_{3n+1,3n+1}=4$,
$l_{3i-1,3i-1}=l_{3i,3i}=4$ for $i\in\{1,2,\ldots,n\}$ and
$l_{3i+1,3i+1}=6$ for $i\in\{1,2,\ldots,n-1\}$.

Assume that $0=\alpha_1< \alpha_2\leq \cdots\leq \alpha_{3n+1}$ and
$0<\beta_1\leq \beta_2\leq \cdots\leq \beta_{3n+1}$ are the roots
of $P_{L_{A}}(x)=0$ and $P_{L_{S}}(x)=0$, respectively.
Note that $$L_{A}\cdot\eta=0,~if~
\eta=(1,1,1,1,1,\ldots,1,1,1,1)'.$$
That is to say, $0$ is one of the eigenvalues of $L_{A}$.  In the line with
(1.1), one gets the following lemma.

\begin{lem} Assume that $O_n$ are the linear crossed octagonal graphs. Then
\begin{eqnarray}
Kf(O_n)&=&2\cdot(3n+1)\cdot\bigg(\sum_{i=2}^{3n+1}\frac{1}{\alpha_i}+\sum_{j=1}^{3n+1}\frac{1}{\beta_j}\bigg).
\end{eqnarray}
\end{lem}

According to the matrix $L_S$, it is straightforward to verified
that
\begin{eqnarray}
\sum_{j=1}^{3n+1}\frac{1}{\beta_j}=(2n+2)\cdot\frac{1}{4}+(n-1)\cdot\frac{1}{6}=\frac{2n+1}{3}.
\end{eqnarray}

In what follows, we will be concerned on the calculations of
$\sum_{i=2}^{3n+1}\frac{1}{\alpha_i}.$ At this place, we firstly give the following two facts which are used to facilitate the proof of Lemma 3.2.
\vskip 0.5 cm

\noindent{\bf Fact 1.} $(-1)^{3n}a_{4n}=(3n+1)\cdot2^{3n}.$
\vskip 0.1 cm

\noindent{\bf Proof.} For the sake of convenience, denoted $W_i$ the $i$-th order
principal submatrix, that is, yield by the first $i$ rows and
columns of $L_{A}, i=1,2,\cdots,3n$. Let $w_i=\det W_i$. Then
$w_1=2$ and $w_2=4$. By expanding $\det W_i$ with regard to its
last row, one obtains $w_i=4w_{i-1}-4w_{i-2}$ for $3\leq i\leq3n$.
By a straight calculation, one gets $w_i=2^i.$

 Notice that $(-1)^{3n}a_{3n}$ is equal to the sum of all
principal minors of $L_{A}$ as discussed above and each of them
both have $3n$ rows and columns, we have
\begin{eqnarray}
(-1)^{3n}a_{3n}&=&\sum_{i=1}^{3n+1}\det
L_{A}[i]=\sum_{i=1}^{3n+1}\det\left(
  \begin{array}{cc}
    W_{i-1} & 0 \\
    0 & U_{3n+1-i} \\
  \end{array}
\right) =\sum_{i=1}^{3n+1}\det W_{i-1} \cdot\det U_{3n+1-i},
\end{eqnarray}

where
\begin{eqnarray*}
 W_{i-1}=
\left(
  \begin{array}{cccc}
    l_{11} & -2 & \cdots & 0\\
    -2 & l_{22} & \cdots & 0 \\
    \vdots & \vdots & \ddots & \vdots \\
    0 & 0 & \cdots  & l_{i-1,i-1} \\
  \end{array}
\right),
\end{eqnarray*}
\begin{eqnarray*}
U_{3n+1-i}=\left(
  \begin{array}{cccc}
    l_{i+1,i+1} & \cdots & 0 & 0\\
    \vdots & \ddots & \vdots & \vdots \\
    0 & \cdots & l_{3n,3n} &-2 \\
    0 & \cdots & -2&l_{3n+1,3n+1} \\
  \end{array}
\right).
\end{eqnarray*}

Let $w_0=1$ and $\det U_0=1$. Assume that the $R[j_1,\ldots,j_k]$
is the submatrix of
$R$, of which is derived by removing the
$j_1$-th, $\ldots$, $j_k$-th rows and columns. According to the
the structure of the symmetry matrix $L_A$, one gets $\det U_{3n+1-i}=\det
W_{3n+1-i}.$ Thus one arrives
\begin{eqnarray*}
(-1)^{3n}a_{3n}&=&\sum_{i=1}^{3n+1}\det
L_{A}[i]=\sum_{i=1}^{3n+1}w_{i-1}w_{3n+1-i}\\
&=&\sum_{i=2}^{3n}w_{i-1}w_{3n+1-i}+2w_{3n}\\
&=&\sum_{i=2}^{3n}2^{i-1}\cdot2^{3n+1-i}+2w_{3n}\\
&=&(3n+1)\cdot2^{3n}.
\end{eqnarray*}
The result as desired.\hfill\rule{1ex}{1ex}
\vskip 0.5 cm

\noindent{\bf Fact 2.}
$(-1)^{3n-1}a_{3n-1}=n(3n+1)(3n+2)\cdot2^{3n-2}.$
\vskip 0.1 cm

\noindent{\bf Proof.} Note that $(-1)^{3n-1}a_{3n-1}$ is equal
to the sum of all principal minors of $L_{A}$ as discussed above
and each of them both have $3n-1$ rows and columns, one
immediately obtains that
\begin{eqnarray*}
L_{A}[i,j]=\left(
  \begin{array}{ccc}
  W_{i-1}& 0&0\\
    0&X_{j-1-i} &0\\
   0& 0& U_{3n+1-j}\\
  \end{array}
\right), 1\leq i<j\leq 3n+1,
\end{eqnarray*}

where
\begin{eqnarray*}
 W_{i-1}=
\left(
  \begin{array}{cccc}
    l_{11} & -2 & \cdots & 0\\
    -2 & l_{22} & \cdots & 0 \\
    \vdots & \vdots & \ddots & \vdots \\
    0 & 0 & \cdots  & l_{i-1,i-1} \\
  \end{array}
\right),~
 X_{j-1-i}=
\left(
  \begin{array}{cccc}
    l_{i+1,i+1} & -2 & \cdots & 0\\
    -2 & l_{i+2,i+2} & \cdots & 0 \\
    \vdots & \vdots & \ddots & \vdots \\
    0 & 0 & \cdots  & l_{j-1,j-1} \\
  \end{array}
\right),
\end{eqnarray*}

\begin{eqnarray*}
 U_{3n+1-j}=
\left(
  \begin{array}{cccc}
    l_{j+1,j+1} & \cdots & 0 & 0\\
    \vdots & \ddots & \vdots & \vdots \\
    0 & \cdots & l_{3n,3n} &-2 \\
    0 & \cdots & -2 &l_{3n+1,3n+1} \\
  \end{array}
\right),
\end{eqnarray*}
where $\det U_{3n+1-j}=1$, if $j=3n+1.$

In view of the matrix $L_A$, one has
\begin{eqnarray*}
X_{j-1-i}=
\left(
  \begin{array}{cccc}
    4 & -2 & \cdots & 0\\
    -2 & 4 & \cdots & 0 \\
    \vdots & \vdots & \ddots & \vdots \\
    0 & 0 & \cdots  & 4 \\
  \end{array}
\right).
\end{eqnarray*}

Let $x_i:=\det X_i.$ We immediately get $x_i=4x_{i-1}-4x_{i-2}$
and $x_1=4.$ Hence, we can obtain $x_i=(1+i)\cdot2^i.$

Note that
\begin{eqnarray*}
(-1)^{3n-1}a_{3n-1}&=&\sum_{1\leq i< j}^{3n+1}\det
L_{A}[i,j]=\sum_{1\leq i< j}^{3n+1}\det W_{i-1}\det X_{j-1-i}\det
U_{3n+1-j}\\
&=&\sum_{1\leq i< j}^{3n+1}w_{i-1}\cdot w_{3n+1-j}\cdot
x_{j-1-i}\\
&=&n(3n+1)(3n+2)\cdot2^{3n-2}.
\end{eqnarray*}

The proof has completed.\hfill\rule{1ex}{1ex}

\begin{lem} Assume that $0=\alpha_1< \alpha_2\leq \cdots\leq \alpha_{3n+1}$
are the eigenvalues of $L_{A}$. One gets
\begin{eqnarray*}
\sum_{i=2}^{3n+1}\frac{1}{\alpha_i}=\frac{n(3n+2)}{4}.
\end{eqnarray*}
\end{lem}

\noindent{\bf Proof.} Notice that
$$P_{L_{A}}(x)=\det(xI-L_{A})=x(x^{3n}+a_{1}\cdot
x^{3n-1}+\cdots+a_{3n-1}\cdot x+a_{3n})$$ with $a_{3n}\neq 0.$
One can easily get that $\frac{1}{\alpha_2}, \frac{1}{\alpha_3},\cdots,
\frac{1}{\alpha_{3n+1}}$ are the roots of the following expression

$$a_{3n}\cdot x^{3n}+a_{3n-1}\cdot x^{3n-1}+\cdots+a_{1}\cdot
x+1=0.$$

In the line with the Vieta's theorem, we have
\begin{eqnarray*}
\sum_{i=2}^{3n+1}\frac{1}{\alpha_i}=\frac{(-1)^{3n-1}\cdot
a_{3n-1}}{(-1)^{3n}\cdot a_{3n}}.
\end{eqnarray*}
Together with Facts 1 and 2, Lemma 3.2
immediately holds.\hfill\rule{1ex}{1ex}


Combining with Lemma 3.2 and (3.6), one gets the desired result of Theorem 1.1.
\hfill\rule{1ex}{1ex}

In what follows,
we list all the Kirchhoff indices of linear crossed
octagonal graph in Table 1, where $1\leq
n\leq 20.$
\begin{table}[h]
\setlength{\abovecaptionskip}{0.05cm} \centering\vspace{.3cm}
\caption{The Kirchhoff indices from $O_{1}$
to $O_{20}$}
\begin{tabular}{c|c|c|c|c|c|c|c|c|c}
  \hline
  $G$ & $Kf(G)$ & $G$ & $Kf(G)$ & $G$ & $Kf(G)$ & $G$ & $Kf(G)$ & $G$ & $Kf(G)$ \\
  \hline
  $O_{1}$ & $18.00$ & $O_{5}$ & $797.33$ & $O_{9}$ & $4008.67$ & $O_{13}$ & $11380.00$ & $O_{17}$ & $24639.33$ \\
  $O_{2}$ & $79.33$ & $O_{6}$ & $1304.67$ & $O_{10}$ & $5394.00$ & $O_{14}$ & $14075.33$ & $O_{18}$ & $29076.67$
  \\
  $O_{3}$ & $211.67$ & $O_{7}$ & $1991.00$ & $O_{11}$ & $7066.33$ & $O_{15}$ & $17165.67$ & $O_{19}$ & $ 34017.00$
  \\
  $O_{4}$ & $442.00$ & $O_{8}$ & $2883.33$ & $O_{12}$ & $9052.67$ & $O_{16}$ & $20678.00$ & $O_{20}$ & $ 39487.33$\\
  \hline
\end{tabular}
\end{table}

In the following, we will devote to prove the ratio between the
Kirchhoff and Wiener index of linear crossed octagonal graph.
\vskip 0.5 cm
\noindent{\bf Proof of Theorem 1.2.} For distance $d_{ij}$, there are four different
species of vertex $i$ while we calculate the Wiener index of linear crossed octagonal graphs. The
expressions of each type of $i$ are listed as below.
\begin{itemize}
\item Vertex 1 of $O_n$:
\begin{eqnarray*}
g_1(n)=1+\sum_{k=1}^{3n}2k=1+3n+9n^2.
\end{eqnarray*}
\item Vertex $3s-1~(s=1,2,\ldots,n)$ of $O_n$:
\begin{eqnarray*}
g_2(n)=\sum_{s=1}^n\bigg(2+\sum_{k=1}^{3s-2}2k+\sum_{k=1}^{3n+2-3s}2k\bigg)=2n(1+3n+3n^2).
\end{eqnarray*}
\item Vertex $3s~(s=1,2,\ldots,n)$ of $O_n$:
\begin{eqnarray*}
g_3(n)=\sum_{s=1}^n\bigg(2+\sum_{k=1}^{3s-1}2k+\sum_{k=1}^{3n+1-3s}2k\bigg)=n(1+3n+6n^2).
\end{eqnarray*}
\item Vertex $3s+1~(s=1,2,\ldots,n-1)$ of $O_n$:
\begin{eqnarray*}
g_4(n)=\sum_{s=1}^{n-1}\bigg(1+\sum_{k=1}^{3s}2k+\sum_{k=1}^{3n-3s}2k\bigg)=-1+n-6n^2+6n^3.
\end{eqnarray*}
\end{itemize}
Thus, one has
\begin{eqnarray*}
W(O_n)=\frac{4g_1(n)+2g_2(n)+2g_3(n)+2g_4(n)}{2}=1+10n+21n^2+18n^3.
\end{eqnarray*}
Combining with $W(O_n)$ and $Kf(O_n)$, one gets
\begin{eqnarray*}
\lim_{n\rightarrow \infty}\frac{Kf(O_n)}{W(O_n)}=\frac{1}{4}.
\end{eqnarray*}
The result as desired.\hfill\rule{1ex}{1ex}

\section{Proofs of Theorems 1.3, 1.4 and 1.5}
\ \ \ \ In this section, we will first give the proof of Theorem 1.3. According to (1.2), one concerns about the normalized Laplacian spectrum of $O_n$ due to its crucial role in the proof of Theorem 1.3. One then considers the number of spanning trees of linear crossed graph in Theorem 1.4. The ratio between multiplicative degree-Kirchhoff and Gutman index is determined in Theorem 1.5.

\noindent{\bf Proof of Theorem 1.3.}
According to (2.4), one has
\begin{eqnarray*}
 \mathcal {L}_A&=&
\left(
  \begin{array}{cccccccccc}
\frac{2}{3} & -\frac{1}{\sqrt{3}} & 0& 0& 0& \cdots& 0& 0& 0& 0\\
-\frac{1}{\sqrt{3}} & 1 & -\frac{1}{2}& 0& 0& \cdots& 0& 0& 0& 0\\
0& -\frac{1}{2} & 1 & -\frac{1}{\sqrt{5}} & 0& \cdots& 0& 0& 0& 0\\
0& 0& -\frac{1}{\sqrt{5}}&\frac{4}{5} &-\frac{1}{\sqrt{5}} &  \cdots& 0& 0& 0& 0\\
0& 0& 0&-\frac{1}{\sqrt{5}} & 1 &\cdots  & 0& 0& 0& 0\\
\vdots& \vdots& \vdots& \vdots& \vdots&\ddots  & \vdots& \vdots& \vdots& \vdots\\
0& 0& 0& 0& 0& \cdots & \frac{4}{5}&-\frac{1}{\sqrt{5}} & 0& 0\\
0& 0& 0& 0& 0& \cdots& -\frac{1}{\sqrt{5}}&1 &-\frac{1}{2} & 0\\
0& 0& 0& 0& 0& \cdots& 0&-\frac{1}{2} &1 &-\frac{1}{\sqrt{3}}  \\
0& 0& 0& 0& 0& \cdots& 0& 0&-\frac{1}{\sqrt{3}} &\frac{2}{3}  \\
  \end{array}
\right)_{(3n+1)\times (3n+1)},
\end{eqnarray*}
where $l_{11}=l_{3n+1,3n+1}=\frac{2}{3}$,
$l_{3i-1,3i-1}=l_{3i,3i}=1$ for $i\in\{1,2,\ldots,n\}$,
$l_{3i+1,3i+1}=\frac{4}{5}$ for $i\in\{1,2,\ldots,n-1\}$. Also,
$l_{12}=l_{21}=l_{3n,3n+1}=l_{3n+1,3n}=-\frac{1}{\sqrt{3}}$,
$l_{2+3i,3+3i}=l_{3+3i,2+3i}=-\frac{1}{2}$, for
$i\in\{0,1,2,\ldots,n-1\}$,
$l_{4+3i,3(i+1)}=l_{3(i+1),4+3i}=l_{5+3i,1+3(i+1)}=l_{1+3(i+1),5+3i}=-\frac{1}{\sqrt{5}}$, for
$k\in\{0,1,2,\ldots,n-1\}$.
\begin{eqnarray*}
 \mathcal{L}_S&=&
\left(
  \begin{array}{cccccccccc}
\frac{4}{3} & 0 & 0& 0& 0& \cdots& 0& 0& 0& 0\\
0 & 1 & 0& 0& 0& \cdots& 0& 0& 0& 0\\
0& 0 & 1 & 0 & 0& \cdots& 0& 0& 0& 0\\
0& 0& 0&\frac{6}{5} &0 &  \cdots& 0& 0& 0& 0\\
0& 0& 0&0 & 1 &\cdots  & 0& 0& 0& 0\\
\vdots& \vdots& \vdots& \vdots& \vdots&\ddots  & \vdots& \vdots& \vdots& \vdots\\
0& 0& 0& 0& 0& \cdots & \frac{6}{5}&0 & 0& 0\\
0& 0& 0& 0& 0& \cdots& 0&1 &0 & 0\\
0& 0& 0& 0& 0& \cdots& 0&0 &1 &0  \\
0& 0& 0& 0& 0& \cdots& 0& 0&0 &\frac{4}{3}  \\
  \end{array}
\right)_{(3n+1)\times (3n+1)},
\end{eqnarray*}
where $l_{11}=l_{3n+1,3n+1}=\frac{4}{3}$,
$l_{3i-1,3i-1}=l_{3i,3i}=1$ for $i\in\{1,2,\ldots,n\}$ and
$l_{3i+1,3i+1}=\frac{6}{5}$ for $i\in\{1,2,\ldots,n-1\}$.

Assume that $0=\lambda_1< \lambda_2\leq \cdots\leq \lambda_{3n+1}$ and
$0<\delta_1\leq \delta_2\leq \cdots\leq \delta_{3n+1}$ are the
roots of $P_{\mathcal{L}_{A}}(x)=0$ and
$P_{\mathcal{L}_{s}}(x)=0$, respectively. Notice that
$$\mathcal{L}_{A}\cdot\eta=0,~if~
\eta=(3,2\sqrt{3},2\sqrt{3},\sqrt{15},\ldots,\sqrt{15},2\sqrt{3},2\sqrt{3},3)'.$$
Namely, $0$ is one of the eigenvalues of $\mathcal{L}_{A}$.  By
\emph{(1.2)}, one get the lemma as below.

\begin{lem} Assume that $O_n$ are the linear crossed octagonal graphs. Then
\begin{eqnarray}
Kf^*(O_n)=2\cdot(13n+1)\cdot(\sum_{i=2}^{3n+1}\frac{1}{\lambda_i}+\sum_{j=1}^{3n+1}\frac{1}{\delta_j}).
\end{eqnarray}
\end{lem}

According to the matrix $\mathcal{L}_S$, one can easily
check that
\begin{eqnarray}
\sum_{j=1}^{3n+1}\frac{1}{\delta_j}=2n+2\cdot\frac{3}{4}+(n-1)\cdot\frac{5}{6}=\frac{17n+4}{6}.
\end{eqnarray}

In what follows, we will be concerned on the calculation of
$\sum_{i=2}^{3n+1}\frac{1}{\lambda_i}.$ At this point, we firstly give the following two facts which are used to facilitate the proof of Lemma 4.2.

\vskip 0.5 cm
\noindent{\bf Fact 3.}
$(-1)^{3n}b_{3n}=\frac{13n+1}{18}\big(\frac{1}{10}\big)^{n-1}.$
\vskip 0.1 cm

\noindent{\bf Proof.} For the sake of convenience, denoted $M_i$ the $i$-th order
principal submatrix, that is, yield by the first $i$ rows and
columns of $\mathcal{L}_{A}, i=1,2,\cdots,3n$.

Let $m_i=\det M_i$. Then
$$m_1=\frac{2}{3}, m_2=\frac{1}{3}, m_3=\frac{1}{6},
 m_4=\frac{1}{15}, m_5=\frac{1}{30}, m_6=\frac{1}{60},$$ and
\begin{eqnarray*}
\begin{cases}
m_{3k}=m_{3k-1}-\frac{1}{4}m_{3k-2};\\
m_{3k+1}=\frac{4}{5}m_{3k}-\frac{1}{5}m_{3k-1};\\
m_{3k+2}=m_{3k+1}-\frac{1}{5}m_{3k}.
\end{cases}
\end{eqnarray*}

By a explicit calculation, one can get those formulas to further
forms as below.
\begin{eqnarray*}
\begin{cases}
m_{3k}=\frac{5}{3}\cdot(\frac{1}{10})^k,~~for~
1\leq k\leq n;\\
m_{3k+1}=\frac{2}{3}\cdot(\frac{1}{10})^k,~~for~0\leq k\leq n-1;\\
m_{3k+2}=\frac{1}{3}\cdot(\frac{1}{10})^k,~~for~0\leq k\leq n-1.
\end{cases}
\end{eqnarray*}

 Note that $(-1)^{3n}b_{3n}$ is equal to the
sum of all principal minors of $\mathcal{L}_{A}$ as discussed
above and each of them both have $3n$ rows and columns, we have
\begin{eqnarray}
(-1)^{3n}b_{3n}=\sum_{i=1}^{3n+1}\det
\mathcal{L}_{A}[i]=\sum_{i=1}^{3n+1}\det\left(
  \begin{array}{cc}
    M_{i-1} & 0 \\
    0 & R_{3n+1-i} \\
  \end{array}
\right) =\sum_{i=1}^{3n+1}\det M_{i-1} \cdot\det R_{3n+1-i},
\end{eqnarray}

where
\begin{eqnarray*}
 M_{i-1}=
\left(
  \begin{array}{cccc}
    l_{11} & -\frac{1}{\sqrt{3}} & \cdots & 0\\
    -\frac{1}{\sqrt{3}} & l_{22} & \cdots & 0 \\
    \vdots & \vdots & \ddots & \vdots \\
    0 & 0 & \cdots  & l_{i-1,i-1} \\
  \end{array}
\right),
\end{eqnarray*}
\begin{eqnarray*}
R_{3n+1-i}=\left(
  \begin{array}{cccc}
    l_{i+1,i+1} & \cdots & 0 & 0\\
    \vdots & \ddots & \vdots & \vdots \\
    0 & \cdots & l_{3n,3n} &-\frac{1}{\sqrt{3}} \\
    0 & \cdots & -\frac{1}{\sqrt{3}}&l_{3n+1,3n+1} \\
  \end{array}
\right).
\end{eqnarray*}

Let $m_0=1$ and $\det R_0=1$. Since a permutation similarity
transformation of a square matrix preserves its determinant and
the structure of matrix $\mathcal{L}_A$, one gets $\det
R_{3n+1-i}=\det M_{3n+1-i}.$ In the line with (4.10), one gets
\begin{eqnarray*}
(-1)^{3n}b_{3n}&=&\sum_{i=1}^{3n+1}\det
\mathcal{L}_{A}[i]=\sum_{i=2}^{3n}\det
\mathcal{L}_{A}[i]+2m_{3n}\\
&=&\sum_{s=1}^{n}\det \mathcal{L}_{A}[3s]+\sum_{s=1}^{n-1}\det
\mathcal{L}_{A}[3s+1]
+\sum_{s=0}^{n-1}\det \mathcal{L}_{A}[3s+2]+2m_{3n}\\
&=&\sum_{s=1}^{n}m_{3(s-1)+2}m_{3(n-s)+1}+\sum_{s=1}^{n-1}m_{3k}m_{3(n-k)}
+\sum_{s=0}^{n-1}m_{3s+1}m_{3(n-s-1)+2}+2m_{3n}\\
&=&\frac{13n+1}{18}\big(\frac{1}{10}\big)^{n-1}.
\end{eqnarray*}
The desired result holds.\hfill\rule{1ex}{1ex}

\vskip 0.5 cm
\noindent{\bf Fact 4.}
$(-1)^{3n-1}b_{3n-1}=\frac{1}{72}(169n^3+39n^2+22n)\big(\frac{1}{10}\big)^{n-1}.$
\vskip 0.1 cm

\noindent{\bf Proof.} Notice that $(-1)^{3n-1}b_{3n-1}$ is equal
to the sum of all principal minors of $\mathcal{L}_{A}$ as
discussed above and each of them both have $3n-1$ rows and
columns, that is
\begin{eqnarray*}
\mathcal{L}_{A}[i,j]=\left(
  \begin{array}{ccc}
  M_{i-1}& 0&0\\
    0&Y_{j-1-i} &0\\
   0& 0& R_{3n+1-j}\\
  \end{array}
\right), 1\leq i<j\leq 3n+1,
\end{eqnarray*}

where
 \begin{eqnarray*}
 M_{i-1}=
\left(
  \begin{array}{cccc}
    c_{11} & -\frac{1}{\sqrt{3}} & \cdots & 0\\
    -\frac{1}{\sqrt{3}} & c_{22} & \cdots & 0 \\
    \vdots & \vdots & \ddots & \vdots \\
    0 & 0 & \cdots  & c_{i-1,i-1} \\
  \end{array}
\right),~
 Y_{j-1-i}=
\left(
  \begin{array}{cccc}
    c_{i+1,i+1} & c_{i+1,i+2} & \cdots & 0\\
    c_{i+2,i+1} & c_{i+2,i+2} & \cdots & 0 \\
    \vdots & \vdots & \ddots & \vdots \\
    0 & 0 & \cdots  & c_{j-1,j-1} \\
  \end{array}
\right),
\end{eqnarray*}

\begin{eqnarray*}
 R_{3n+1-j}=
\left(
  \begin{array}{cccc}
    c_{j+1,j+1} & \cdots & 0 & 0\\
    \vdots & \ddots & \vdots & \vdots \\
    0 & \cdots & c_{3n,3n} &-\frac{1}{\sqrt{3}} \\
    0 & \cdots & -\frac{1}{\sqrt{3}} &c_{3n+1,3n+1} \\
  \end{array}
\right),
\end{eqnarray*}
where $\det U_{3n+1-j}=1$, if $j=3n+1.$

Note that
\begin{eqnarray}
(-1)^{3n-1}b_{3n-1}=\sum_{1\leq i< j}^{3n+1}\det
\mathcal{L}_{A}[i,j]=\sum_{1\leq i< j}^{3n+1}m_{i-1}m_{3n+1-j}\det
Y_{j-1-i}.
\end{eqnarray}

According to the (4.11), all these cases of $i$ and $j$ are listed as follows.

{\bf{Case 1.}} If $i\equiv0~(mod~3), j\equiv0~(mod~3)$, $1\leq i<
j\leq 3n+1$. That is, $i=3p,~j=3q$ and $1\leq p< q\leq n$. In this
case, $Y$ is a square matrix with order $3q-3p-1$.
\begin{eqnarray*}
 \det Y&=&
\left|
  \begin{array}{ccccccc}
    \frac{4}{5} & -\frac{1}{\sqrt{5}} & 0 & 0 & \cdots & 0 & 0\\
    -\frac{1}{\sqrt{5}} & 1 & -\frac{1}{2}& 0 & \cdots & 0 & 0\\
     0 & -\frac{1}{2} & 1 &-\frac{1}{\sqrt{5}} & \cdots & 0 & 0\\
    0 & 0 & -\frac{1}{\sqrt{5}} & \frac{4}{5} & \cdots & 0 & 0\\
    \vdots & \vdots & \vdots & \vdots & \ddots & \vdots & \vdots\\
    0 & 0 & 0 & 0 & \cdots & \frac{4}{5} & -\frac{1}{\sqrt{5}}\\
    0 & 0 & 0 & 0 & \cdots & -\frac{1}{\sqrt{5}} & 1   \\
  \end{array}
\right|=\frac{3(q-p)}{5}\cdot\bigg(\frac{1}{10}\bigg)^{q-p-1}.
\end{eqnarray*}

{\bf{Case 2.}} If $i\equiv0~(mod~3), j\equiv1~(mod~3)$, $1\leq i<
j\leq 3n+1$. That is, $i=3p,~j=3q+1$ and $1\leq p\leq q\leq n-1$.
In this case, $Y$ is a square matrix with order $3q-3p$.
\begin{eqnarray*}
 \det Y&=&
\left|
  \begin{array}{ccccccc}
    \frac{4}{5} & -\frac{1}{\sqrt{5}} & 0 & 0 & \cdots & 0 & 0\\
    -\frac{1}{\sqrt{5}} & 1 & -\frac{1}{2}& 0 & \cdots & 0 & 0\\
     0 & -\frac{1}{2} & 1 &-\frac{1}{\sqrt{5}} & \cdots & 0 & 0\\
    0 & 0 & -\frac{1}{\sqrt{5}} & \frac{4}{5} & \cdots & 0 & 0\\
    \vdots & \vdots & \vdots & \vdots & \ddots & \vdots & \vdots\\
    0 & 0 & 0 & 0 & \cdots & 1 & -\frac{1}{2}\\
    0 & 0 & 0 & 0 & \cdots & -\frac{1}{2} & 1   \\
  \end{array}
\right|=(3q-3p+1)\cdot\bigg(\frac{1}{10}\bigg)^{q-p}.
\end{eqnarray*}

{\bf{Case 3.}} If $i\equiv0~(mod~3), j\equiv2~(mod~3)$, $1\leq i<
j\leq 3n+1$. That is, $i=3p,~j=3q+2$ and $1\leq p\leq q\leq n-1$.
In this case, $Y$ is a square matrix with order $3q-3p+1$.
\begin{eqnarray*}
 \det Y&=&
\left|
  \begin{array}{ccccccc}
    \frac{4}{5} & -\frac{1}{\sqrt{5}} & 0 & 0 & \cdots & 0 & 0\\
    -\frac{1}{\sqrt{5}} & 1 & -\frac{1}{2}& 0 & \cdots & 0 & 0\\
     0 & -\frac{1}{2} & 1 &-\frac{1}{\sqrt{5}} & \cdots & 0 & 0\\
    0 & 0 & -\frac{1}{\sqrt{5}} & \frac{4}{5} & \cdots & 0 & 0\\
    \vdots & \vdots & \vdots & \vdots & \ddots & \vdots & \vdots\\
    0 & 0 & 0 & 0 & \cdots & 1 & -\frac{1}{\sqrt{5}}\\
    0 & 0 & 0 & 0 & \cdots & -\frac{1}{\sqrt{5}} & \frac{4}{5}   \\
  \end{array}
\right|=(\frac{6}{5}q-\frac{6}{5}p+\frac{4}{5})\cdot\bigg(\frac{1}{10}\bigg)^{q-p}.
\end{eqnarray*}

{\bf{Case 4.}} If $i\equiv0~(mod~3), j=3n+1$, $1\leq i\leq 3n$.
That is, $i=3p$ and $1\leq p\leq n$. In this case, $Y$ is a square
matrix with order $3q-3p$.
\begin{eqnarray*}
 \det Y&=&
\left|
  \begin{array}{ccccccc}
    \frac{4}{5} & -\frac{1}{\sqrt{5}} & 0 & 0 & \cdots & 0 & 0\\
    -\frac{1}{\sqrt{5}} & 1 & -\frac{1}{2}& 0 & \cdots & 0 & 0\\
     0 & -\frac{1}{2} & 1 &-\frac{1}{\sqrt{5}} & \cdots & 0 & 0\\
    0 & 0 & -\frac{1}{\sqrt{5}} & \frac{4}{5} & \cdots & 0 & 0\\
    \vdots & \vdots & \vdots & \vdots & \ddots & \vdots & \vdots\\
    0 & 0 & 0 & 0 & \cdots & 1 & -\frac{1}{2}\\
    0 & 0 & 0 & 0 & \cdots & -\frac{1}{2} & 1   \\
  \end{array}
\right|=(3n-3p+1)\cdot\bigg(\frac{1}{10}\bigg)^{n-p}.
\end{eqnarray*}

{\bf{Case 5.}} If $i\equiv1~(mod~3), j\equiv0~(mod~3)$, $1\leq i<
j\leq 3n+1$. That is, $i=3p+1,~j=3q$ and $0\leq p< q\leq n$. In
this case, $Y$ is a square matrix with order $3q-3p-2$.
\begin{eqnarray*}
 \det Y&=&
\left|
  \begin{array}{ccccccc}
    1 & -\frac{1}{2} & 0 & 0 & \cdots & 0 & 0\\
    -\frac{1}{2} & 1 & -\frac{1}{\sqrt{5}}& 0 & \cdots & 0 & 0\\
     0 & -\frac{1}{\sqrt{5}} & \frac{4}{5} &-\frac{1}{\sqrt{5}} & \cdots & 0 & 0\\
    0 & 0 & -\frac{1}{\sqrt{5}} & 1 & \cdots & 0 & 0\\
    \vdots & \vdots & \vdots & \vdots & \ddots & \vdots & \vdots\\
    0 & 0 & 0 & 0 & \cdots & \frac{4}{5} & -\frac{1}{\sqrt{5}}\\
    0 & 0 & 0 & 0 & \cdots & -\frac{1}{\sqrt{5}} & 1   \\
  \end{array}
\right|=(\frac{3}{2}q-\frac{3}{2}p-\frac{1}{2})\cdot\bigg(\frac{1}{10}\bigg)^{q-p-1}.
\end{eqnarray*}

{\bf{Case 6.}} If $i\equiv1~(mod~3), j\equiv1~(mod~3)$, $1\leq i<
j\leq 3n+1$. That is, $i=3p+1,~j=3q+1$ and $0\leq p< q\leq n-1$. In
this case, $Y$ is a square matrix with order $3q-3p-1$.
\begin{eqnarray*}
 \det Y&=&
\left|
  \begin{array}{ccccccc}
    1 & -\frac{1}{2} & 0 & 0 & \cdots & 0 & 0\\
    -\frac{1}{2} & 1 & -\frac{1}{\sqrt{5}}& 0 & \cdots & 0 & 0\\
     0 & -\frac{1}{\sqrt{5}} & \frac{4}{5} &-\frac{1}{\sqrt{5}} & \cdots & 0 & 0\\
    0 & 0 & -\frac{1}{\sqrt{5}} & 1 & \cdots & 0 & 0\\
    \vdots & \vdots & \vdots & \vdots & \ddots & \vdots & \vdots\\
    0 & 0 & 0 & 0 & \cdots & 1 & -\frac{1}{2}\\
    0 & 0 & 0 & 0 & \cdots & -\frac{1}{2} & 1   \\
  \end{array}
\right|=\frac{3}{4}(q-p)\cdot\bigg(\frac{1}{10}\bigg)^{q-p-1}.
\end{eqnarray*}

{\bf{Case 7.}} If $i\equiv1~(mod~3), j\equiv2~(mod~3)$, $1\leq i<
j\leq 3n+1$. That is, $i=3p+1,~j=3q+2$ and $0\leq p\leq q\leq n-1$.
In this case, $Y$ is a square matrix with order $3q-3p$.
\begin{eqnarray*}
 \det Y&=&
\left|
  \begin{array}{ccccccc}
    1 & -\frac{1}{2} & 0 & 0 & \cdots & 0 & 0\\
    -\frac{1}{2} & 1 & -\frac{1}{\sqrt{5}}& 0 & \cdots & 0 & 0\\
     0 & -\frac{1}{\sqrt{5}} & \frac{4}{5} &-\frac{1}{\sqrt{5}} & \cdots & 0 & 0\\
    0 & 0 & -\frac{1}{\sqrt{5}} & 1 & \cdots & 0 & 0\\
    \vdots & \vdots & \vdots & \vdots & \ddots & \vdots & \vdots\\
    0 & 0 & 0 & 0 & \cdots & 1 & -\frac{1}{\sqrt{5}}\\
    0 & 0 & 0 & 0 & \cdots & -\frac{1}{\sqrt{5}} & \frac{4}{5}   \\
  \end{array}
\right|=(3q-3p+1)\cdot\bigg(\frac{1}{10}\bigg)^{q-p}.
\end{eqnarray*}

{\bf{Case 8.}} If $i\equiv1~(mod~3), j=3n+1$, $1\leq i\leq 3n$.
That is, $i=3p+1$ and $0\leq p\leq n-1$. In this case, $Y$ is a
square matrix with order $3n-3p-1$.
\begin{eqnarray*}
 \det Y&=&
\left|
  \begin{array}{ccccccc}
    1 & -\frac{1}{2} & 0 & 0 & \cdots & 0 & 0\\
    -\frac{1}{2} & 1 & -\frac{1}{\sqrt{5}}& 0 & \cdots & 0 & 0\\
     0 & -\frac{1}{\sqrt{5}} & \frac{4}{5} &-\frac{1}{\sqrt{5}} & \cdots & 0 & 0\\
    0 & 0 & -\frac{1}{\sqrt{5}} & 1 & \cdots & 0 & 0\\
    \vdots & \vdots & \vdots & \vdots & \ddots & \vdots & \vdots\\
    0 & 0 & 0 & 0 & \cdots & 1 & -\frac{1}{2}\\
    0 & 0 & 0 & 0 & \cdots & -\frac{1}{2} & 1   \\
  \end{array}
\right|=\frac{3}{4}(n-p)\cdot\bigg(\frac{1}{10}\bigg)^{n-p-1}.
\end{eqnarray*}

{\bf{Case 9.}} If $i\equiv2~(mod~3), j\equiv0~(mod~3)$, $1\leq i<
j\leq 3n+1$. That is, $i=3p+2,~j=3q$ and $0\leq p< q\leq n$. In
this case, $Y$ is a square matrix with order $3q-3p-3$.
\begin{eqnarray*}
 \det Y&=&
\left|
  \begin{array}{ccccccc}
    1 & -\frac{1}{\sqrt{5}} & 0 & 0 & \cdots & 0 & 0\\
    -\frac{1}{\sqrt{5}} & \frac{4}{5} & -\frac{1}{\sqrt{5}}& 0 & \cdots & 0 & 0\\
     0 & -\frac{1}{\sqrt{5}} & 1 &-\frac{1}{2} & \cdots & 0 & 0\\
    0 & 0 & -\frac{1}{2} & 1 & \cdots & 0 & 0\\
    \vdots & \vdots & \vdots & \vdots & \ddots & \vdots & \vdots\\
    0 & 0 & 0 & 0 & \cdots & \frac{4}{5} & -\frac{1}{\sqrt{5}}\\
    0 & 0 & 0 & 0 & \cdots & -\frac{1}{\sqrt{5}} & 1   \\
  \end{array}
\right|=(3q-3p-2)\cdot\bigg(\frac{1}{10}\bigg)^{q-p-1}.
\end{eqnarray*}

{\bf{Case 10.}} If $i\equiv2~(mod~3), j\equiv1~(mod~3)$, $1\leq i<
j\leq 3n+1$. That is, $i=3p+2,~j=3q+1$ and $0\leq p< q\leq n-1$. In
this case, $Y$ is a square matrix with order $3q-3p-2$.
\begin{eqnarray*}
 \det Y&=&
\left|
  \begin{array}{ccccccc}
    1 & -\frac{1}{\sqrt{5}} & 0 & 0 & \cdots & 0 & 0\\
    -\frac{1}{\sqrt{5}} & \frac{4}{5} & -\frac{1}{\sqrt{5}}& 0 & \cdots & 0 & 0\\
     0 & -\frac{1}{\sqrt{5}} & 1 &-\frac{1}{2} & \cdots & 0 & 0\\
    0 & 0 & -\frac{1}{2} & 1 & \cdots & 0 & 0\\
    \vdots & \vdots & \vdots & \vdots & \ddots & \vdots & \vdots\\
    0 & 0 & 0 & 0 & \cdots & 1 & -\frac{1}{2}\\
    0 & 0 & 0 & 0 & \cdots & -\frac{1}{2} & 1   \\
  \end{array}
\right|=(\frac{3}{2}q-\frac{3}{2}p-\frac{1}{2})\cdot\bigg(\frac{1}{10}\bigg)^{q-p-1}.
\end{eqnarray*}

{\bf{Case 11.}} If $i\equiv2~(mod~3), j\equiv2~(mod~3)$, $1\leq i<
j\leq 3n+1$. That is, $i=3p+2,~j=3q+2$ and $0\leq p< q\leq n-1$. In
this case, $Y$ is a square matrix with order $3q-3p-1$.
\begin{eqnarray*}
 \det Y&=&
\left|
  \begin{array}{ccccccc}
    1 & -\frac{1}{\sqrt{5}} & 0 & 0 & \cdots & 0 & 0\\
    -\frac{1}{\sqrt{5}} & \frac{4}{5} & -\frac{1}{\sqrt{5}}& 0 & \cdots & 0 & 0\\
     0 & -\frac{1}{\sqrt{5}} & 1 &-\frac{1}{2} & \cdots & 0 & 0\\
    0 & 0 & -\frac{1}{2} & 1 & \cdots & 0 & 0\\
    \vdots & \vdots & \vdots & \vdots & \ddots & \vdots & \vdots\\
    0 & 0 & 0 & 0 & \cdots &  1& -\frac{1}{\sqrt{5}}\\
    0 & 0 & 0 & 0 & \cdots & -\frac{1}{\sqrt{5}} & \frac{4}{5}   \\
  \end{array}
\right|=(\frac{3}{5}q-\frac{3}{5}p)\cdot\bigg(\frac{1}{10}\bigg)^{q-p-1}.
\end{eqnarray*}

{\bf{Case 12.}} If $i\equiv2~(mod~3), j=3n+1$, $1\leq i\leq 3n$.
That is, $i=3p+2$ and $0\leq p\leq n-1$. In this case, $Y$ is a
square matrix with order $3n-3p-2$.
\begin{eqnarray*}
 \det Y&=&
\left|
  \begin{array}{ccccccc}
    1 & -\frac{1}{\sqrt{5}} & 0 & 0 & \cdots & 0 & 0\\
    -\frac{1}{\sqrt{5}} & \frac{4}{5} & -\frac{1}{\sqrt{5}}& 0 & \cdots & 0 & 0\\
     0 & -\frac{1}{\sqrt{5}} & 1 &-\frac{1}{2} & \cdots & 0 & 0\\
    0 & 0 & -\frac{1}{2} & 1 & \cdots & 0 & 0\\
    \vdots & \vdots & \vdots & \vdots & \ddots & \vdots & \vdots\\
    0 & 0 & 0 & 0 & \cdots & 1 & -\frac{1}{2}\\
    0 & 0 & 0 & 0 & \cdots & -\frac{1}{2} & 1   \\
  \end{array}
\right|=(\frac{3}{2}n-\frac{3}{2}p-\frac{1}{2})\cdot\bigg(\frac{1}{10}\bigg)^{n-p-1}.
\end{eqnarray*}

Thus, one has
\begin{eqnarray}
(-1)^{3n-1}b_{3n-1}&=&\sum_{1\leq i\leq
j}^{3n+1}m_{i-1}m_{3n+1-j}\cdot\det Y_{j-1-i} \nonumber \\
&=&E_1+E_2+E_3,
\end{eqnarray}

where
\begin{eqnarray*}
E_1&=&\sum_{1\leq p<q\leq n}\det\mathcal{L}_{A}[3p,3q]
+\sum_{1\leq p\leq q\leq n-1}\det\mathcal{L}_{A}[3p,3q+1]
+\sum_{1\leq p\leq q\leq n-1}\det\mathcal{L}_{A}[3p,3q+2]\\
&&+\sum_{1\leq p\leq n}\det\mathcal{L}_{A}[3p,3n+1]\\
&=&\frac{n(n-1)(n+1)}{45}\bigg(\frac{1}{10}\bigg)^{n-2}+\frac{5n(n-1)^2}{18}\bigg(\frac{1}{10}\bigg)^{n-1}+
\frac{n^2(n-1)}{45}\bigg(\frac{1}{10}\bigg)^{n-2}
+\frac{n(3n-1)}{6}\bigg(\frac{1}{10}\bigg)^{n-1}\\
&=&\frac{13n^3-5n^2-2n}{18}\bigg(\frac{1}{10}\bigg)^{n-1}.
\end{eqnarray*}

\begin{eqnarray*}
E_2&=&\sum_{1\leq p<q\leq
n}\det\mathcal{L}_{A}[3p+1,3q]+\sum_{1\leq p< q\leq
n-1}\det\mathcal{L}_{A}[3p+1,3q+1] +\sum_{1\leq p\leq q\leq
n-1}\det\mathcal{L}_{A}[3p+1,3q+2]\\
&&+\sum_{1\leq p\leq
n-1}\det\mathcal{L}_{A}[3p+1,3n+1]+\sum_{1\leq q\leq
n}\det\mathcal{L}_{A}[1,3q]+\sum_{1\leq q\leq
n-1}\det\mathcal{L}_{A}[1,3q+1]\\
&&+\sum_{0\leq
q\leq n-1}\det\mathcal{L}_{A}[1,3q+2]+X_{3n-1}\\
&=&\frac{5n^2(n-1)}{18}\bigg(\frac{1}{10}\bigg)^{n-1}+\frac{25n(n-1)(n-2)}{72}\bigg(\frac{1}{10}\bigg)^{n-1}
+\frac{5n(n-1)^2}{18}\bigg(\frac{1}{10}\bigg)^{n-1}
+\frac{5n(n-1)}{8}\bigg(\frac{1}{10}\bigg)^{n-1}\\
&&+\frac{n(3n+1)}{6}\bigg(\frac{1}{10}\bigg)^{n-1}+
\frac{5n(n-1)}{8}\bigg(\frac{1}{10}\bigg)^{n-1}+\frac{n(3n-1)}{6}\bigg(\frac{1}{10}\bigg)^{n-1}+\frac{3n}{4}\bigg(\frac{1}{10}\bigg)^{n-1}\\
&=&\frac{65n^3+27n^2+34n}{72}\bigg(\frac{1}{10}\bigg)^{n-1}.
\end{eqnarray*}

\begin{eqnarray*}
E_3&=&\sum_{0\leq p<q\leq
n}\det\mathcal{L}_{A}[3p+2,3q]+\sum_{0\leq p< q\leq
n-1}\det\mathcal{L}_{A}[3p+2,3q+1]+\sum_{0\leq p< q\leq
n-1}\det\mathcal{L}_{A}[3p+2,3q+2]\\
&&+\sum_{0\leq p\leq n-1}\det\mathcal{L}_{A}[3p+2,3n+1]\\
&=&\frac{2n^2(n+1)}{9}\bigg(\frac{1}{10}\bigg)^{n-1}+\frac{5n^2(n-1)}{18}\bigg(\frac{1}{10}\bigg)^{n-1}+
\frac{n(n-1)(n+1)}{45}\bigg(\frac{1}{10}\bigg)^{n-2}+\frac{n(3n+1)}{6}\bigg(\frac{1}{10}\bigg)^{n-1}\\
&=&\frac{13n^3+8n^2-n}{18}\bigg(\frac{1}{10}\bigg)^{n-1}.
\end{eqnarray*}

Substituting $E_1, E_2$ and $E_3$ to (4.12), we obtain
\begin{eqnarray*} (-1)^{3n-1}b_{3n-1}&=&
\frac{169n^3+39n^2+22n}{72}\bigg(\frac{1}{10}\bigg)^{n-1}.
\end{eqnarray*}
The proof of Fact 4. completed.\hfill\rule{1ex}{1ex}

\begin{lem} Assume that $0=\lambda_1< \lambda_2\leq \cdots\leq \lambda_{3n+1}$
are the eigenvalues of $\mathcal{L}_{A}$. One has
\begin{eqnarray}
\sum_{i=2}^{3n+1}\frac{1}{\lambda_i}=\frac{169n^3+39n^2+22n}{4\cdot(13n+1)}.
\end{eqnarray}
\end{lem}

\noindent{\bf Proof.}
Note that
$$P_{\mathcal{L}_{A}}(x)=\det(xI-\mathcal{L}_{A})=x(x^{3n}+a_{1}\cdot
x^{3n-1}+\cdots+a_{3n-1}\cdot x+a_{3n})$$ with $a_{3n}\neq
0.$ One gets that $\frac{1}{\alpha_1}, \frac{1}{\alpha_3},\cdots,
\frac{1}{\alpha_{3n+1}}$ satisfies the following expression

$$a_{3n}\cdot x^{3n}+a_{3n-1}\cdot x^{3n-1}+\cdots+a_{1}\cdot
x+1=0.$$

In the line with the Vieta's theorem, we have
\begin{eqnarray*}
\sum_{i=2}^{3n+1}\frac{1}{\lambda_i}=\frac{(-1)^{3n-1}\cdot
b_{3n-1}}{(-1)^{3n}\cdot b_{3n}}.
\end{eqnarray*}
Combining with Facts 3 and 4, the proof of Lemma 4.2 has completed. \hfill\rule{1ex}{1ex}

Substituting (4.9) and (4.13) to (4.8), we can get
the following theorems.\hfill\rule{1ex}{1ex}


In what follows, we list the multiplicative degree-Kirchhoff indices of linear crossed
octagonal graphs in Table 2, where $1\leq
n\leq 20.$
\begin{table}[h]
\setlength{\abovecaptionskip}{0.05cm} \centering\vspace{.3cm}
\caption{The degree-Kirchhoff indices from $O_{1}$
to $O_{12}$}
\begin{tabular}{c|c|c|c|c|c|c|c|c|c}
  \hline
  $G$ & $Kf^*(G)$ & $G$ & $Kf^*(G)$ & $G$ & $Kf^*(G)$ & $G$ & $Kf^*(G)$ & $G$ & $Kf^*(G)$ \\
  \hline
  $O_{1}$ & $213.00$ & $O_{5}$ & $13063.00$ & $O_{9}$ & $69454.33$ & $O_{13}$ & $201835.00$ & $O_{17}$ & $442653.00$ \\
  $O_{2}$ & $1118.00$ & $O_{6}$ & $21811.33$ & $O_{10}$ & $94158.00$ & $O_{14}$ & $250606.00$ & $O_{18}$ & $523603.33$
  \\
  $O_{3}$ & $3223.33$ & $O_{7}$ & $33788.00$ & $O_{11}$ & $124118.00$ & $O_{15}$ & $306661.33$ & $O_{19}$ & $ 613866.00$
  \\
  $O_{4}$ & $7036.00$ & $O_{8}$ & $49500.00$ & $O_{12}$ & $159841.33$ & $O_{16}$ & $370508.00$ & $O_{20}$ & $ 713948.00$\\
  \hline
\end{tabular}
\end{table}

\noindent{\bf Proof of Theorem 1.4.} Based on the Lemma 2.2, one has $
\prod_{i=1}^{|V(O_n)|}d_i\prod_{i=2}^{3n+1}\alpha_i\prod_{j=1}^{3n+1}\beta_j=2|E(O_n)|\tau(O_n).
$ Then
\begin{eqnarray*}
\prod_{i=1}^{|V(O_n)|}d_i&=&3^4\cdot4^{4n}\cdot5^{2n-2},\\
%
\prod_{i=2}^{3n+1}\alpha_i&=&(-1)^{3n}a_{3n}=\frac{13n+1}{18}\bigg(\frac{1}{10}\bigg)^{n-1},\\
%
\prod_{j=1}^{4n+1}\beta_j&=&\bigg(\frac{4}{3}\bigg)^2\bigg(\frac{6}{5}\bigg)^{n-1} .
\end{eqnarray*}
Hence \begin{eqnarray*}
\tau(O_n)&=&2^{8n+2}\cdot3^{n-1}.
\end{eqnarray*}
The result immediately holds.\hfill\rule{1ex}{1ex}

At this point, all the number of spanning trees of linear
crossed octagonal graphs are listed in Table 3, where
$1\leq n\leq 8.$
\begin{table}[htbp]
\setlength{\abovecaptionskip}{0.05cm}
 \centering \vspace{.3cm}
\caption{All the number of spanning trees from $O_{1}$ to
$O_{8}$} \begin{tabular}{c|c|c|c|c|c|c|c}
  \hline
  $G$ & $\tau(G)$ & $G$ & $\tau(G)$ & $G$ & $\tau(G)$ & $G$ & $\tau(G)$ \\
  \hline
  $O_{1}$ & $1024$ & $O_{3}$ & $603979776$ & $O_{5}$ & $356241767399424$ & $O_{7}$ & $210119944214597861376$  \\
  $O_{2}$ & $786432$ & $O_{4}$ & $463856467968$ & $O_{6}$ & $273593677362757632$ & $O_{8}$ & $161372117156811157536768$
  \\
  \hline
\end{tabular}
\end{table}

\noindent{\bf Proof of Theorem 1.5.} For distance $d_{ij}$, there are four different
species of vertex $i$ while we calculate the Gutman index of linear crossed octagonal graphs. The
expressions of each type of $i$ are listed as follows.
\begin{itemize}
\item Vertex 1 of $O_n$:
\begin{eqnarray*}
f_1(n)&=&\sum_{k\in I}3\cdot4\cdot k\cdot2+\sum_{k\in J}3\cdot5\cdot k\cdot2
+\sum_{k\in I}3\cdot4\cdot (k+1)\cdot2+3\cdot3\cdot(1+3n+3n)\\&=&9\cdot(1+n+13n^2),
\end{eqnarray*}
where $I=\{1,4,\ldots,3n-2\}$ and $J=\{3,6,\ldots,3n-3\}.$
\item Vertex $3s-1~(s=1,2,\ldots,n)$ of $O_n$:
\begin{eqnarray*}
f_2(n)&=&\sum_{s=1}^n\sum_{k=1}^s4\cdot4\cdot 3k\cdot2+\sum_{s=1}^{n}\sum_{k=1}^{n-s}4\cdot4\cdot 3k\cdot2+
\sum_{s=1}^n\sum_{k=1}^{s-1}4\cdot4\cdot (3k-1)\cdot2\\&&+\sum_{s=1}^{n}\sum_{k=1}^{n-s+1}4\cdot4\cdot (3k-2)\cdot2
+\sum_{s=1}^n\sum_{k=1}^{s-1}4\cdot5\cdot (3k-2)\cdot2+\sum_{s=1}^{n}\sum_{k=1}^{n-s}4\cdot5\cdot (3k-1)\cdot2\\
&&+\sum_{s=1}^n4\cdot3(3s-1-1+3n+1-3s+1)\cdot2+4\cdot4\cdot 2\\
&=&4\cdot(8+9n+15n^2+26n^3).
\end{eqnarray*}
\item Vertex $3s~(s=1,2,\ldots,n)$ of $O_n$:
\begin{eqnarray*}
f_3(n)&=&\sum_{s=1}^n\sum_{k=1}^{s-1}4\cdot4\cdot 3k\cdot2+\sum_{s=1}^{n}\sum_{k=1}^{n-s}4\cdot4\cdot 3k\cdot2+
\sum_{s=1}^n\sum_{k=1}^{s}4\cdot4\cdot (3k-2)\cdot2\\&&+\sum_{s=1}^{n}\sum_{k=1}^{n-s}4\cdot4\cdot (3k-1)\cdot2
+\sum_{s=1}^n\sum_{k=1}^{s-1}4\cdot5\cdot (3k-1)\cdot2+\sum_{s=1}^{n}\sum_{k=1}^{n-s}4\cdot5\cdot (3k-2)\cdot2\\
&&+\sum_{s=1}^n4\cdot3(3s-1+3n+1-3s)\cdot2+4\cdot4\cdot 2\\
&=&4\cdot(8-3n+3n^2+26n^3).
\end{eqnarray*}
\item Vertex $3s-2~(s=2,3,\ldots,n)$ of $O_n$:
\begin{eqnarray*}
f_4(n)&=&\sum_{s=2}^n\sum_{k=1}^{s-2}5\cdot5\cdot 3k\cdot2+\sum_{s=2}^{n}\sum_{k=1}^{n-s}5\cdot5\cdot 3k\cdot2+
\sum_{s=2}^n\sum_{k=1}^{s-1}5\cdot4\cdot (3k-1)\cdot2\\&&+\sum_{s=2}^{n}\sum_{k=1}^{n-s+1}5\cdot4\cdot (3k-2)\cdot2
+\sum_{s=2}^n\sum_{k=1}^{s-1}4\cdot5\cdot (3k-2)\cdot2+\sum_{s=2}^{n}\sum_{k=1}^{n-s+1}4\cdot5\cdot (3k-1)\cdot2\\
&&+\sum_{s=2}^n5\cdot3(3s-2-1+3n+1-3s+2)\cdot2+5\cdot5\\
&=&5\cdot(5+10n-36n^2+26n^3).
\end{eqnarray*}
\end{itemize}
Thus, one has
\begin{eqnarray*}
Gut(O_n)&=&\frac{4f_1(n)+2f_2(n)+2f_3(n)+2f_4(n)}{2}\\&=&107+92n+126n^2+338n^3.
\end{eqnarray*}
Combining with $Gut(O_n)$ and $Kf^*(O_n)$, one gets
\begin{eqnarray*}
\lim_{n\rightarrow \infty}\frac{Kf^*(O_n)}{Gut(O_n)}=\frac{1}{4}.
\end{eqnarray*}
The result as desired.\hfill\rule{1ex}{1ex}

\section{Conclusion} \ \ \ \ \ In this paper, we first considered the Laplacian and normalized Laplacian of $O_n.$ Those have some great effects on computing the Kirchhoff and multiplicative degree-Kirchhoff index. The main results, namely, those resistance distance-based graph invariants and the number of spanning trees of $O_n$ were determined. For extended study, one may focus on the computing the resistance distance between any two vertices of $O_n$ on the basis of resistance sum rules or the Kirchhoff, multiplicative degree-Kirchhoff index and the number of spanning trees of linear octagonal M$\ddot{o}$bius graph.

\section*{Acknowledgments}

\ \ \ \ \ This work is supported by National Natural Science
Foundation of China (Nos. 11601006, 11801007).



\end{document}